\documentclass[multphys,vecphys]{svmult}

\usepackage{makeidx}         
\usepackage{graphicx}        
\usepackage{multicol}        
\usepackage[bottom]{footmisc}

\makeindex             

\begin{document}

\title*{Poisson Statistics for the Largest Eigenvalues in Random Matrix Ensembles}
\titlerunning{Poisson Statistics in Random Matrices} 
\author{Alexander Soshnikov
}
\institute{University of California at Davis\\
Department of Mathematics\\
Davis, CA 95616, USA\\
\texttt{soshniko@math.ucdavis.edu}}
%
%
\maketitle


\section{Introduction}
\label{sec:1}


The two archetypal  ensembles of random matrices are Wigner real symmetric (Hermitian) random matrices and Wishart sample covariance real 
(complex) random matrices. In this paper we study the statistical properties of the largest eigenvalues of such matrices in the case when 
the second moments of matrix entries are infinite.
In the first two subsections we consider Wigner ensemble of random matrices and its generalization --
band random matrices.

\subsection{Wigner Random Matrices}
\label{sec:2}

A real symmetric Wigner random matrix is defined as a square symmetric
$ n\times n $ matrix with i.i.d. entries up from the diagonal
\begin{equation}
\label{matrica}
 A=(a_{jk}), \ \ a_{jk}=a_{kj}, \ 1 \leq j \leq k \leq n,  \ \ \{a_{jk}\}_{j<k} - \ {\rm i.i.d.} \ \ {\rm real} \ \ 
{\rm random} \ \ {\rm variables}.
\end{equation}
The diagonal entries $\{a_{ii}\}, \  1\leq i \leq n, $ are usually assumed to be i.i.d. random variables, independent from the off-diagonal 
entries.
A Hermitian Wigner random matrix is defined in a similar way, namely as a square $ n \times n $ Hermitian matrix with i.i.d. entries up
from the diagonal
\begin{equation}
\label{hermitian}
 A=(a_{jk}), \ \ a_{jk}=\overline{a_{kj}}, \ 1 \leq j \leq k \leq n, 
\ \ \{a_{jk}\}_{j< k} - \ {\rm i.i.d.} \ \ {\rm complex} \ \ {\rm random} \ \ {\rm variables}.
\end{equation}
As in the real symmetric case, it is  usually assumed that the diagonal entries 
$ \{ a_{ii} \}, \  1\leq i \leq n, $ are i.i.d. (real) random variables independent
from the off-diagonal entries.

Ensembles (\ref{matrica}) and (\ref{hermitian}) were introduced in mathematical physics by Eugene Wigner in the 1950s (\cite{Wig1}, 
\cite{Wig2}, \cite{Wig3}). Wigner viewed these ensembles as  a  mathematical model to study the statistics of the excited energy levels 
of heavy nuclei.  

The famous Wigner's semicircle law can be formulated as follows. Let the matrix entries in (\ref{matrica}) 
or (\ref{hermitian}) be centered random variables with the tail of distribution decaying sufficiently fast, so that all moments exist.
Denote by $ \lambda_1 \geq \lambda_2 \geq \ldots \geq \lambda_n $ the eigenvalues of a random matrix $n^{-1/2}\*A.$ 
Then the empirical distribution function of the eigenvalues converges, as $n\to \infty$, to a non-random limit

\begin{equation}
\label{semicircle}
\frac{1}{n} \* \# (\lambda_i \leq x, \ \ 1\leq i \leq n)  \to F(x)= \int_{-\infty}^x f(t) \* dt,
\end{equation}
where the density of the semicircle law is given by $ f(t)= \frac{1}{\pi \* \sigma^2} \* \sqrt{2 \sigma^2 - x^2}, \ $ 
for $ t \in [-\sqrt{2}\* \sigma, \sqrt{2}\*\sigma], $ and $\sigma^2 $ is the second moment of matrix entries.

This result was subsequently strengthened by many mathematicians (see e.g. \cite{Ar}, \cite{Pas}, \cite{FK}). In its general form 
(due to Pastur and Girko), the theorem holds if the matrix entries  of $A$  satisfy the Lindeberg-Feller condition:
$ \ \ \frac{1}{n^2} \sum_{1\leq i \leq j \leq n} \int_{|x|>\tau \* \sqrt{n}} x^2\* dF_{ij}(x) \to 0, \ \ $  where $ F_{ij}(x) $ is the 
distribution function
of $ a_{ij}^{(n)}.$ 

From the analytical point of view, the simplest examples of  Wigner random matrices are given by the so-called 
Gaussian Orthogonal and Unitary Ensembles (GOE and GUE for short).  The GUE is defined as the ensemble of  $n\times n $ Hermitian matrices
with the Gaussian entries $\Re a_{jk} \sim N(0, 1/2), \ \ \Im a_{jk} \sim N(0, 1/2), \ \ 1 \leq j < k \leq n; \ \ a_{ii}\sim N(0, 1), \ \ 
1\leq i \leq n$ (see (\cite{M}, chapter 6). The joint distribution of the matrix entries has the form
\begin{equation}
\label{GUE}
P(dA)= const_n \* \exp\left(-\frac{1}{2}\*Tr(A^2)\right) \* dA, 
\end{equation}
where $dA= \prod_{j\leq k} d \Re a_{ij} \* d \Im a_{jk} \* \prod_{i=1}^n da_{ii}\ $ is the Lebesgue measure 
on the space of $n$-dimensional Hermitian matrices. The joint distribution of the eigenvalues is given by its density
\begin{equation}
\label{gue}
p_n(x_1, \ldots, x_n)= Z_n^{-1}\* \prod_{1\leq i < j \leq n} (x_i-x_j)^2 \* \exp\left(-\frac{1}{2} \* \sum_i x_i^2\right).
\end{equation}

The normalization constants in (\ref{GUE}) and (\ref{gue}) are known. What is more, one can calculate explicitely the
$k$-point correlation functions (see \cite{M}, chapter 6). This allows one to study the local distribution of the eigenvalues, both in the 
bulk of ths spectrum and at its edges in great detail. In particular, a celebrated result of Tracy and Widom (see \cite{TW1}) states that
\begin{equation}
\label{Tracy2}
\lim_{n \to \infty} \Pr \left( \lambda_{max} \leq 2\*\sqrt{n}  + \frac{s}{n^{1/6}} \right) = F_{2}(s)=
\exp\left( - \int_s^{+\infty}(x-s)\* q^2(x)\* dx \right),
\end{equation}
where $ q(x) \ $ is the solution  of the Painl\'{e}ve II differential equation $$ q''(x)=x\*q(x) +2\*q^3(x) $$ with the asymptotics
at infinity$ \ q(x) \sim Ai(x) \ $ at $ x = +\infty.$

The limiting $k$-point correlation function at the edge of spectrum is given by the formula
\begin{equation}
\label{corrgue}
\rho_k(x_1,\ldots, x_k)= \det \left( K(x_i,x_j) \right)_{1 \leq i,j \leq k},
\end{equation}
where
\begin{equation}
\label{airyker}
K(x,y)=K_{Airy}(x,y) = \frac{Ai(x)\*Ai'(y)-Ai'(x)\*Ai(y)}{x-y}
\end{equation}
is a so-called Airy kernel. We refer the reader to  \cite{TW1} and \cite{For}  for the details.
We recall that the $k$-point correlation function is defined in such a way that
for any disjoint subintervals of the real line $ \ I_1, I_2, \ldots, I_k, \ $ one has
$$ \ E \prod_{i=1}^k  \#(I_i) = \int_{I_1} \ldots \int_{I_k} \rho_k(x_1, \ldots, x_k) \* dx_1 \ldots dx_k, $$
where $ \ \#(I) \ $ denotes the number of the eigenvalues in $ I. \ $ A probabilistic interpretation of the above formula is that
$ \ \ \rho_k(x_1, \ldots, x_k) \*dx_1 \ldots dx_k \ $ is the probability to find an eigenvalue in each of the $k$ infinitesimal intervals
$ \ [x_i, x_i +dx_i], \ i=1, \ldots, k.$

The Gaussian Orthogonal Ensemble (GOE) is defined as the ensemble of $n\times n$ Wigner real symmetric 
random matrices with the Gaussian entries. More precisely, we assume that $ a_{ij}, \ \ 1\leq i \leq j \leq n, \ $ are 
independent  Gaussian $\ N(0, 1+\delta_{ij}) $
random variables (see e.g. \cite{M}, chapter 7). The joint distribution of the matrix entries has the form
\begin{equation}
\label{GOE}
P(dA)= c_n \* \exp\left(-\frac{1}{4}\* Tr(A^2)\right) \* dA, 
\end{equation}
where $dA= \prod_{i\leq j} da_{ij} \ $ is the Lebesgue measure 
on the space of $n$-dimensional real symmetric matrices.
The distribution (\ref{GOE}) induces the joint distribution of the eigenvalues of the GOE matrix, given by its density

\begin{equation}
\label{goe}
p_n(x_1, \ldots, x_n)= Z_n^{-1}\* \prod_{1\leq i < j \leq n} |x_i-x_j| \* \exp\left(-\frac{1}{4} \* \sum_i x_i^2\right).
\end{equation}

The limiting distribution of the (normalized) largest eigenvalue of a GOE matrix was calculated by Tracy and Widom in (\cite{TW2}).
\begin{equation}
\label{Tracy}
\lim_{n \to \infty} \Pr \left( \lambda_{max} \leq 2\*\sqrt{n} + \frac{s}{n^{1/6}} \right) = F_{1}(s)=
\exp\left( -\frac{1}{2} \* \int_s^{+\infty} q(x) + (x-s)\* q^2(x)\* dx \right),
\end{equation}
The Tracy-Widom distribution (\ref{Tracy}) was obtained by studying the asymptotic properties of the $k$-point correlation functions at 
the edge of the spectrum. 
The $k-$point correlation function in the GOE ensemble has the pfaffian form. In the limit $ n \to \infty$ the $k$-point correlation function
at the edge of the spectrum  is given by  the following formula

\begin{equation}
\label{corrgoe}
\rho_k(x_1,\ldots, x_k)= \left(\det \left( K(x_i,x_j) \right)_{1 \leq i,j \leq k}\right)^{1/2},
\end{equation}
where $K(x,y)$ is a $2\times 2$ matrix kernel such that
\begin{eqnarray}
K_{11}(x,y)&=&K_{22}(y,x)= K_{Airy}(x,y) +\frac{1}{2} \* Ai(x) \*\int_{-\infty}^y Ai(t)\* dt, \\
K_{12}(x,y)&=&-\frac{1}{2} \* Ai(x)\* Ai(y) - \frac{\partial}{\partial y} K_{Airy}(x,y),  \\
K_{21}(x,y)&=& \int_0^{+\infty} \left( \int_{x+u}^{+\infty} Ai(v)\*dv\right) \* Ai(x+u) \* du - \epsilon(x-y) +\frac{1}{2} 
\int_y^x Ai(u) \* du  \nonumber \\
&+& \frac{1}{2} \*\int_x^{+\infty} Ai(u) \* du \* \int_{-\infty}^y Ai(v) \* dv,
\end{eqnarray} 
where $ \epsilon(z)= \frac{1}{2} \* sign(z). \ $

\subsection{Band Random Matrices}
A band random matrix is a generalization of a Wigner random matrix ensemble (\ref{matrica}), (\ref{hermitian}).
A real symmetric (aperiodic) band random matrix is defined as a square symmetric
$ n\times n $ matrix $A=(a_{jk}) \ $ such that $a_{ij}=0 \ $ unless $|i-j| \leq d_n, \ $  and
\begin{equation}
\label{matrica1}
\{a_{jk}, \ j\leq k; \ |j-k| \leq d_n \} - \ {\rm i.i.d.} \ \ {\rm real} \ \ 
{\rm random} \ \ {\rm variables}.
\end{equation}

A Hermitian band random matrix is defined in a similar way, namely as a square $ n \times n $ Hermitian matrix $A=(a_{jk}), \ $ such that
$a_{ij}=0 \ $ unless $|i-j| \leq d_n, \ $ and
\begin{equation}
\label{hermitian1}
\{a_{jk}, \ j \leq k; \ |j-k| \leq d_n \} - \ {\rm i.i.d.} \ \ {\rm complex} \ \ 
{\rm random} \ \ {\rm variables}.
\end{equation}

If $d_n=n-1,$ we obtain the Wigner ensemble of random matrices. 
A matrix is called a periodic band matrix if $|i-j| $ is replaced above by $|i-j|_1=\min( |i-j|, n-|i-j|). \  $
Band random matrices have been studied in the last fifteen years (see for example \cite{CMI}, \cite{CG}, \cite{MPK}, \cite{Gui}). 
In the periodic case, the limiting distribution of the eigenvalues
of $d_n^{-1}\*A$ is given by the semi-circle law, provided matrix entries have a finite second moment. In the aperiodic case, the limiting 
distribution of the eigenvalues is different from the semi-circle law, unless $d_n/n \to 0 $ (see e.g. \cite{MPK}). One of the most 
interesting problems involving band random matrices is the localization/ delocalization properties of the eigenvalues.  It is conjectured in 
physical literature, that the eigenvalues of band random matrices are localized if $d_n =O(n^{1/2}).\ $ As far as we know, there are no 
rigorous results yet in this direction.

\subsection{Sample Covariance Random Matrices}
\label{sec:3}

Sample covariance random matrices have been studied in mathematical statistics for the last seventy-five years. We refer to \cite{Mui},
\cite{Wil} and \cite{Ja} for the applications of  spectral properties of Wishart random matrices
in multivariate statistical analysis.

Let $A$ be a large $ m \times n $ real rectangular random matrix with  independent  identically distributed entries. In applications, one 
is often interested in
the statistical behavior of the singular values of $A$ in the limit $m\to \infty, \ n \to \infty.$  This is equivalent to studying the 
eigenvalues of a positive-definite matrix $M=A^t\*A\ $ in the limit of large dimensions. 
Without loss of generality, one can assume that $m \geq n \ $(since the spectrum of $A\*A^t \ $ differs from the spectrum of $A^t\*A \ $
only by a zero eigenvalue of multiplicity $m-n.$

The analogue of the Wigner semicircle law was proved by Marchenko and Pastur (\cite{MP}).  Let $m \to \infty, \ n \to \infty \ $ in such a way
that $m/n \to \gamma \geq 1. \ $ Assume  $ E |a_{ij}|^{2+\epsilon} < +\infty, \ $ where $ \ \epsilon>0 \ $ is an arbitrary positive number. 
Then the empirical distribution function of the eigenvalues
of $\frac{1}{m}\* A^t\*A$ converges to a non-random limit, known as the Marchenko-Pastur distribution

\begin{equation}
\label{pastur}
\frac{1}{n} \* \#( \lambda_i \leq x, \ i=1, \ldots,n) \to G_{\gamma}(x)= \int_{-\infty}^x g_{\gamma}(t) \* dt,
\end{equation}
where the spectral density $g(t)$ is supported on the interval $[a,b], \ \ a= \sigma^2\*(1-\gamma^{-1/2})^2, 
\ b=\sigma^2\*(1+\gamma^{-1/2})^2,  \ \sigma^2=E a_{11}^2, \ $ and
$g(t)= \frac{1}{2\*\pi\*t\*\gamma\*\sigma^2}\sqrt{(b-t)(t-a)}, \ \ t\in [a,b].$

The case $a_{ij} \sim N(0,1) \ \ 1\leq i,j \leq n, \ $ is known in the literature as the Wishart (Laguerre) ensemble of real 
sample covariance matrices. The joint distribution of the eigenvalues of $M$ is defined by its density.
Similarly to the Gaussian ensembles of real symmetric and Hermitian matrices discussed in Subsection 1.1, many important statistical 
quantities 
in the Wishart ensemble can be calculated explicitely. 
For example, the joint probability density of the  eigenvalues is given by the formula
\begin{equation}
\label{realwishart}
p_n(x_1, \ldots, x_n)= Z_{n,m}^{-1}\* \prod_{1\leq i < j \leq n} |x_i-x_j| \* \prod_{i=1}^n x_i^{m-n-1} \*
\exp(-x_i/2).
\end{equation}
It was shown by Johnstone (\cite{J}), that the largest eigenvalue of a Wishart random matrix converges, after a proper rescaling, to 
the Tracy-Widom distribution $F_1.$ Namely, let $m \to \infty, \ n \to \infty, \ m/n \to \gamma$ and
$\mu_{m,n}= (n^{1/2}+m^{1/2})^2, \ \ \sigma_{m,n}=(n^{1/2}+m^{1/2})\*(n^{-1/2}+m^{-1/2})^{1/3}.\ $ Then

\begin{equation}
\label{john}
\Pr \left( \lambda_{max}(A^t\*A) \leq \mu_{m,n} + s \* \sigma_{m,n} \right) \to F_1(s)
\end{equation}

One can also show (see \cite{So2}), that the rescaled $k$-point correlation function at the edge of the spectrum converge in the limit 
to (\ref{corrgoe}).

Finally, we want to remark, that there is a long-standing interest in nuclear physics in the spectral properties of the complex sample 
covariance matrices $A^*\*A,$ where the entries of a reactangular matrix $A$ are independent identically distributed complex random variables
(see e.g. \cite{Wig3}, \cite{Br}, \cite{FS}, \cite{qcd}, \cite{Been}). We refer the reader to \cite{So2} and the references therein for 
additional information.

\subsection{Universality in Random Matrices}
\label{sec:5}

The universality conjecture in Random Matrix Theory states, loosely speaking, that the local statistical properties of a few eigenvalues 
in the bulk or at the edge of the spectrum are independent of the distribution of individual matrix entries in the limit of large dimension.
The only thing that should matter is, whether the matrix is real symmetric, Hermitian or self-dual quaternion Hermitian.

For Wigner random matrices, the conjecture was rigorously proven at the edge of the spectrum, both for real symmetric and Hermitian case 
in \cite{So1}, provided that all moments of matrix entries exist and do not grow faster than the moments of a 
Gaussian distribution, and the odd moments vanish.  In particular, it was shown that the largest eigenvalue, after proper rescaling, converges
in distribution to the Tracy-Widom law.
In the bulk of the spectrum, the conjecture was proven by Johansson (\cite{Jo1}) for 
Wigner Hermitian matrices, provided the marginal distribution of a matrix entry has a Gaussian component. We refer to \cite{D}
and references therein for the universality results in the  unitary ensembles of random matrices.

The situation for sample covariance random matrices is quite similar (see papers by Soshnikov \cite{So2} and Ben Arous and P\'{e}ch\'{e} 
\cite{BAP}). 

The natural question is  how general such results are? What happens if matrix entries have only a finite number of moments?
In this article we consider the extreme case when the entries of $A$ do not have a finite second moment.
In the next section, we discuss  spectral properties of Wigner random matrices and, more generally, band random matrices
when marginal distribution  of matrix 
entries has heavy tails. As was shown in \cite{So3}, the statistics of the largest eigenvalues of such matrices are given by a Poisson 
inhomogeneous random point process. In Section 3 we discuss a similar result (although in a weaker form) 
for the largest eigenvalues of sample 
covariance random matrices with Cauchy entries. Section 4 is devoted to conclusions.

\section{Wigner and Band Random Matrices with Heavy Tails of Marginal Distributions}
\label{sec:6}

In this section we consider ensembles of Wigner real symmetric and Hermitian matrices (\ref{matrica}) and (\ref{hermitian}), and band 
real symmetric and Hermitian random matrices (\ref{matrica1}), (\ref{hermitian}) with the 
additional condition on the tail of the marginal distribution

\begin{equation}
\label{tail}
G(x)=\Pr (|a_{jk}| > x) = \frac{h(x)}{x^{\alpha}}, 
\end{equation}
where
$ 0 < \alpha < 2 $ and $ h(x) $ is a slowly varying function at infinity in a sense of Karamata (\cite{Kar}, \cite{Sen}). In 
other words, $h(x)$ is a positive function for all $ x>0,$ such that $ \lim_{x \to \infty} \frac{h(t\*x)}{h(x)}=1 $ for all $ t>0.$
The condition (\ref{tail}) means that the distribution of $\ |a_{ij}| \ $ belongs to the domain of the attraction of a stable distribution
with the index $\alpha$ (see e.g.  \cite{IL}, Theorem 2.6.1).

Without loss of generality, we restrict our attention to the real symmetric case. The results in the Hermitian case are practically the same.
Wigner random matrices (\ref{matrica}), (\ref{hermitian}) with the heavy tails (\ref{tail}), in the special case when  limit 
$\lim_{x\to \infty} h(x) >0 \ $ exists, were considered on a physical level of rigor by Cizeau and Bouchaud
in \cite{CB}. They argued, that the typical eigenvalues of $A$ are of the order of $n^{1/\alpha}.$  Cizeau and Bouchaud also suggested a 
formula  for the limiting spectral density of the empirical distribution function of the eigenvalues of
$n^{-1/\alpha}\*A.$  Unlike the Wigner semicircle and Marchenko-Pastur laws, 
the conjectured limiting spectral density is supported on the whole real line. It was given as
\begin{equation}
\label{plotnost}
f(x)= L_{\alpha/2}^{C(x), \beta(x)}(x),
\end{equation}
where $L_{\alpha}^{C,\beta} $ is a density of a centered L\'{e}vy stable distribution  defined through its Fourier 
transform $ \hat{L}(k) :$
\begin{eqnarray}
\label{Levy}
& & L_{\alpha}^{C,\beta}= \frac{1}{2\*\pi}\*\int dk \* \hat{L}(k) \* e^{i\*k\*x}, \\
& & \ln \hat{L}(k)= -C\*|k|^{\alpha}\*\bigl(1 +i\beta\*sgn(k)\*\tan(\pi\*\alpha/2)\bigr),
\end{eqnarray}
and functions $ C(x), \ \ \beta(x) \ $ satisfy a system of integral equations
\begin{eqnarray}
\label{Cbeta}
& & C(x)=\int_{-\infty}^{+\infty} |y|^{\frac{\alpha}{2}-2}\* L_{\alpha/2}^{C(y), \beta(y)}\bigl(x-\frac{1}{y}\bigr) \* dy, \\
& & \beta(x)= \int_{x}^{+\infty}  L_{\alpha/2}^{C(y), \beta(y)}\bigl(x-\frac{1}{y}\bigr)\* dy.
\end{eqnarray}

We would like to draw the reader's attention to the fact that the density in  (\ref{plotnost}) is not a density of a L\'{e}vy stable 
distribution,
since $C(x), \ \ \beta(x) \ $ are  
functions of $x.$   Cizeau and Bouchaud argued, that the density $f(x) \ $ should decay as $  \frac{1}{x^{1+\alpha}} \ $ at
infinity,  thus suggesting that
the largest eigenvalues of $A$ (in the case  $ \ h(x)=const) \ $ should be of order $n^{\frac{2}{\alpha}}, \ $ and not
$ n^{\frac{1}{\alpha}}, \ $  which is the order of  typical eigenvalues. 

Even though originally proven in \cite{So3} in the Wigner case,
the theorem written below holds in the general case of band random real symmetric (or Hermitian) random matrices
(\ref{matrica1}), (\ref{hermitian1}). 

Let $N_n$ be the number of independent (i.e. $i\leq j$), non-zero matrix entries $a_{ij} $
in $A. \ $ In other words, let $N_n= \#( 1\leq i \leq j \leq n, \ \ |i-j|\leq d_n) \ $ in the aperiodic band case, and $N_n=
\#( 1\leq i \leq j \leq n, \ \ |i-j|_1\leq d_n) \ $ in the periodic band case.
It is not difficult to see, 
that $N_n= \frac{n\*(n+1)}{2}\ $ in the Wigner case, $ \ N_n=n\*(d_n+1)$ in the periodic band case, and $\ N_n=n\*(d_n+1)- 
\frac{d_n\*(d_n+1)}{2} 
\ $ in the aperiodic band case.
Let us define a normalization constant $b_n$ in such a way that

\begin{equation}
\label{bn}
\lim_{n \to \infty} \ N_n\* G(b_n\*x) = \frac{1}{x^{\alpha}},
\end{equation}
for all positive $x>0, \ $ where the tail distribution $G$ has been defined in (\ref{tail}). Normalization $b_n$ naturally appears
(see \cite{LLR} and Remark 1 below), when one studies the extremal values of a sequence of $N_n$ independent identically distributed random 
variables (\ref{tail}). In particular, one
can choose
\begin{equation}
\label{b}
b_n= \inf\{ t: G(t-0) \geq \frac{1}{N_n} \geq G(t+0).
\end{equation}
It follows from (\ref{bn}) and (\ref{b}), that  $ \ N_n^{\alpha -\delta} \ll b_n \ll N_n^{\alpha +\delta} \ $
for arbitrary small positive $\delta, \ $ and $ \frac{N_n\* h(b_n)}{b_n^{\alpha}} \to 1 $ as $ n \to \infty.$

Theorem 1 claims that the largest eigenvalues of $A$ have Poisson statistics in the limit $n \to \infty.$

\begin{theorem}

Let $A$ be a band real symmetric (\ref{matrica1}) or Hermitian (\ref{hermitian1})
random matrix  with a heavy tail of the distribution of matrix entries (\ref{tail}).
Then the random point configuration composed of the positive eigenvalues of $b_n^{-1} \* A$  converges in distribution on the cylinder sets
to the inhomogeneous Poisson random point process on $ (0, +\infty)  $ with the intensity 
$ \rho(x)= \frac{\alpha}{x^{1+\alpha}}.$
\end{theorem}

In other words, let $ 0 < x_1 <y_1 <x_2 <y_2 <\ldots x_k <y_k \leq +\infty, $ and $ I_j=(x_j, y_j), \ \ j=1, \ldots k, \ $ be disjoint 
intervals on the 
positive half-line. Then the counting random variables 
$ \#(I_j)= \# ( 1 \leq i \leq n : \lambda_i \in I_j), \ \ j=1, \ldots, k, $ 
are independent in the limit
$ n \to \infty, $ and have a joint Poisson distribution with the parameters $ \mu_j= \int_{I_j} \rho(x) dx, $  i.e.
\begin{equation}
\label{jointP}
\lim_{n\to \infty} \Pr \left(\#(I_j)=s_j, \ j=1, \ldots, k\right)= \prod_{j=1}^k \frac{\mu_j^{s_j}}{s_j!} \* e^{-\mu_j}.
\end{equation}
For the additional information on Poisson random point processes we refer the reader to \cite{DVJ}.

\begin{corollary}
Let  $\lambda_k$ be the $k$-th largest eigenvalue of $b_n^{-1}\*A, $ then
\begin{equation}
\label{kth}
\lim_{n \to \infty} \Pr (\lambda_k \leq x) =  \exp(- \* x^{-\alpha}) \* \sum_{l=0}^{k-1} \frac{x^{-l\*\alpha}}{l!}.
\end{equation}
In particular, $ \lim_{n \to \infty} \Pr (\lambda_1 \leq x) =  \exp(- \* x^{-\alpha}).$

\end{corollary}

{\bf Remark 1}
{\it The equivalent formulation of the theorem is the following. Let $k$ be a finite positive integer. Then the joint distribution of the 
first 
$k$ largest eigenvalues of $b_n^{-1}\*A$ is asymptotically (in the limit $n \to \infty$) the same as the joint distribution of the 
first $k$ order statistics of $ \{ b_n^{-1} \*|a_{ij}|, \ \ 1\leq i \leq j \leq n \} $. It is a classical result,
that extremal values of the sequence of independent identically distributed random variables with heavy tails distributions (\ref{tail})
have Poisson statistics (see e.g. \cite{LLR}, Theorem 2.3.1).}

Theorem 1 was proven in \cite{So3} in the Wigner (i.e. full matrix) case (\ref{matrica}), (\ref{hermitian}). The proof in the general 
(band matrix) case is essentially the same.
However, it should be noted, that the original proof of Theorem 1 in \cite{So3} contained a little 
mistake, which could be easily corrected.
The corrections are
due in two places.

First of all, 
the correct formulation of the part c) of 
Lemma 4 from \cite{So3} (p. 87)  should  state, that
for any positive constant $\delta >0, $ with probability going to 1  there is no no row $ 1 \leq i \leq n, $  that contains at least two 
entries 
greater in absolute value than $b_n^{\frac{3}{4} +\delta}\ $. In other words, the exponent $\frac{1}{2} +\delta \ $ in
$b_n^{\frac{1}{2} +\delta} \ $ in part c) of Lemma 4
must be replaced by
$\frac{3}{4} +\delta.\ $ After this correction, the statement is true. Indeed, the probability that there is a row with at least two entries 
greater than $b_n^{\frac{3}{4} +\delta} \ $
can be estimated from above by $n^3\*\left(G(b_n^{\frac{3}{4} +\delta})\right)^2. \ $ It
follows from (\ref{tail}), (\ref{bn}) and (\ref{b}), that this probability goes to zero.

Also, the formula (28) in Lemma 5 (p. 88) should read
\begin{equation}
\label{netu}
\Pr \{ \exists i, \ \ 1 \leq i \leq n : \ \ \max_{1\leq j \leq n} |a_{ij}| > b_n^{\frac{3}{4}+\frac{\alpha}{8}}, 
\ \ ( \sum_{1 \leq j \leq n}|a_{ij}| ) - 
\max_{1\leq j \leq n} |a_{ij}| >b_n^{\frac{3}{4}+\frac{\alpha}{8}} \} \to 0
\end{equation}
as $ n \to \infty. \ $ In other words, the exponent $ \frac{1}{2}+\frac{\alpha}{4} \ $ in $ \ b_n^{\frac{1}{2}+\frac{\alpha}{4}} \ $
must be replaced by $ \frac{3}{4}+\frac{\alpha}{8}.$ 
The key step of the proof of Lemma 5  was to show, that for any fixed row $i$ and arbitrary small positive $\delta,$
the probability $ \Pr ( \sum_{j: |a_{ij}| \leq b_n^{\frac{1}{2}+\delta}} |a_{ij}| \geq b_n^{\frac{1}{2}+2\*\delta} )\ $ can be estimated 
from above by  $ \exp(-n^{\epsilon}), $ where $\epsilon= \epsilon(\delta, \alpha) >0. \ $  We then concluded, that with probability
going to 1, there is no row $i$ such that $ \sum_{j: |a_{ij}| \leq b_n^{\frac{1}{2}+\delta}} |a_{ij}| \geq b_n^{\frac{1}{2}+2\*\delta}. $ 
To establish 
(\ref{netu}), it is  enough to prove that for any fixed row $i$
\begin{equation}
\label{raz} 
\Pr ( \sum_{j: b_n^{1/2+\delta} \leq |a_{ij}| \leq b_n^{\frac{3}{4}+\delta}} |a_{ij}| \geq b_n^{\frac{3}{4}+2\*\delta} ) < 
\exp(-n^{\epsilon}),
\end{equation}
for sufficiently small positive $\epsilon.$  The proof is very similar to the argument presented in Lemma 5, and is left to the reader.

\section{Real Sample Covariance Matrices with Cauchy Entries}

Let $A$ be a rectangular $m \times n$ matrix with independent identically distributed  entries with the marginal probability 
distribution of matrix entries satisfying (\ref{tail}).
Based on the results in the last section, one can expect that the largest eigenvalues have Poisson statistics as well.
At this point, we have been able to prove it only in a weak form, and only when matrix entries have Cauchy distribution.

We recall, that the probability density of the Cauchy distribution is given by the formula $  f(x)= \frac{1}{\pi (1+x^2)}.$
Cauchy distribution is very important in probability theory (see e.g. \cite{Fel}). In particular, Cauchy distribition is a (1,1,0) stable 
distribution, i.e. the scale parameter is 1, the index of the distribution $\alpha=1$ and the symmetry parameter is zero (see \cite{IL},
\cite{LLR}).

The following theorem was proven by Fyodorov and Soshnikov in \cite{SF}

\begin{theorem}
 Let  $A $ be a random rectangular $m \times n$ matrix ($ m \geq n $ ) with i.i.d. Cauchy entries and let
$ z $ be a complex number with a positive real part.
Then, as $ n \to \infty $
we have
\begin{equation}
\label{glavform}
\lim_{n \to \infty} E \left(\det(1 +\frac{z}{m^2 \*n^2} \* A^t \* A )\right)^{-1/2}  =
\exp\left(-\frac{2}{\pi} \* \sqrt{z} \right)=
{\bf E} \prod_{i=1}^{\infty} (1+ z\* x_i)^{-1/2},
\end{equation}
where
we consider the branch of $\sqrt{z}$ on $D=\{ z: \Re z >0 \} $ such that $ \sqrt{1}=1,$
$ E $ denotes the mathematical expectation with respect to the random matrix ensemble defined above, ${\bf E}$ denotes the mathematical 
expectation with respect to the inhomogeneous Poisson random point process on the positive half-axis with the intensity 
$ \frac{1}{\pi \* x^{3/2}}, $ 
and the convergence 
is uniform inside $D$ (i.e. it is unform on compact subsets of $D$). For  a real positive $z=t^2, \ t \in R^1, $ one can  estimate 
the rate of convergence, namely
\begin{equation}
\label{glavformA}
\lim_{n \to \infty} E \left(\det(1 +\frac{t^2}{m^2 \*n^2} \* A^t \* A )\right)^{-1/2} =
\exp\left(-\frac{2}{\pi} \* |t| \* \big(1+o(n^{-1/2+\epsilon})\big)\right),
\end{equation}
where $ \epsilon $ is an arbitrary small positive number and the convergence is uniform on the compact subsets of
$ [0, +\infty).$

\end{theorem}

The result of Theorem 2 allows a  generalization to the case of 
a sparse random matrix with Cauchy entries.  
Let, as before, $ \{ a_{jk}\}, \ 1 \leq j \leq m, \ 1 \leq k \leq n, \ $ be i.i.d. Cauchy random 
variables, and $ B= (b_{jk}) $ be a $ m \times n $ non-random rectangular $ 0-1$ matrix such that the number of non-zero entries in each 
column
is fixed and equals to $d_n.$  Let $d_n $ grow polynomially, i.e.
$ \  \ b_n \geq n^{\alpha}, $ for some $ 0<\alpha \leq 1. $  Also assume that $ \ \ln m \ $ grows much slower than than any power of $n $.

Define a $ m \times n $ rectangular matrix $A$ with the entries
$ A_{jk}= b_{jk} \* a_{jk}, \ \ 1\leq j \leq m, \ 1 \leq k \leq n. \ $ Let
$ \lambda_1 \geq \lambda_2 \ldots \geq \lambda_n $ denote
the eigenvalues of $ A^t\* A$. The appropriate rescaling for the largest eigenvalues in this case is 
$ \tilde{\lambda_i} = \frac{\lambda_i}{m^2 \* d_n^2}, \ \ i=1, \ldots, n.$

\begin{theorem}
 Let  $A$ be a sparse random rectangular $m \times n$ matrix ($ m \geq n $) defined as above, and let
$ z $ be a complex number with a positive real part.
Then, as $ n \to \infty $
we have
\begin{eqnarray}
\label{glavform1}
& &
\lim_{n \to \infty} E \left(\det(1 +\frac{z}{m^2 \* d_n^2} \* 
A^t \* A )\right)^{-1/2} =
\lim_{n \to \infty} E \prod_{i=1}^n (1+ z\* \tilde{\lambda_i})^{-1/2} \\
&=&
\exp\left(-\frac{2}{\pi} \* \sqrt{z}  \right) 
= {\bf E} \prod_{i=1}^{\infty} (1+ z\* x_i)^{-1/2},
\end{eqnarray}
where, as in Theorem 1.1,
we consider the branch of $\sqrt{z}$ on $D=\{ z: \Re z >0 \} $ such that $ \sqrt{1}=1.$
${\bf E}$ denotes the mathematical 
expectation with respect to the inhomogeneous Poisson random point process on the positive half-axis with the intensity 
$ \frac{1}{\pi \* x^{3/2}},$
and the convergence 
is uniform inside $D$ (i.e. it is unform on the compact subsets of $D$). 
For  a real positive $z=t^2, \ t \in  R^1, $ one can get an estimate 
on the rate of convergence, namely
\begin{equation}
\label{glavform1A}
E \left(\det(1 +\frac{t^2}{m^2 \* d_n^2} \* A^t \* A )\right)^{-1/2} =
\exp\left(-\frac{2}{\pi} \* t \* \big(1+o(d_n^{-1/2+\epsilon})\big)\right),
\end{equation}
where
$ \epsilon $ is an arbitrary small positive number and the convergence is 
uniform on the compact subsets of $[0, +\infty).$

\end{theorem}

The proof relies on the following property of the Gaussian integrals:
\begin{equation}
\label{Gauss}
\bigl(\det(B)\bigr)^{-1/2}=
\left(\frac{1}{\pi}\right)^N \* \int  x \* \exp \bigl(-x \* B \* x^t \bigr)\*d^{N}, 
\end{equation}
where $B$ is an $N$-dimensional matrix with a positive definite Hermitian part (i.e. all eigenvalues of $ B+B^* \ $ are positive),
$x=(x_1, \ldots, x_N) \in R^N, \ $ and $ x \* B \* x^t= \sum_{ij} b_{ij}\*x_i \*x_j.$

Let $B= B(t)= \left( \begin{array}{cc} Id & t i \* A\\ t i \* A^t & Id \end{array} \right).$
Then, one
can write
\begin{equation}
\label{determ}
\bigl(\det(1 + t^2\* A^tA)\bigr)^{-1/2} =
\left( \det \left( \begin{array}{cc} 1 &  ti \* A\\ t i \* A^t & 1 \end{array} \right) \right)^{-1/2}= \bigl(\det(B)\bigr)^{-1/2},
\end{equation}
and apply (\ref{Gauss}) to the r.h.s. of (\ref{determ}). Assuming that the entries of $A$ are independent, one can significantly simplify
the expression, using the fact that entries of $A$ appear linearly in $B(t)\ $(see Proposition 1  of \cite{SF}).  In the Cauchy case, one 
can simplify the calculations even further, and prove that
$ \ \lim_{n \to \infty} E \left(\det(1 +\frac{z}{m^2 \* d_n^2} \* 
A^t \* A )\right)^{-1/2} \ $
exists and equals
$ \exp\left(-\frac{2}{\pi} \* \sqrt{z}  \right). \ $

On the other side, for Poisson random point processes the mathematical expectations of the type  $ {\bf E} \prod_{i=1}^{\infty} (1+f(x_i))
\ $ can be calculated explicitely
\begin{eqnarray}
\label{poisson}
& & {\bf E} \prod_{i=1}^{\infty} (1+f(x_i)) = 1 +\sum_{k=1}^{\infty} {\bf E} \sum_{1 \leq i_1 < i_2 <\ldots < i_k} 
\prod_{j=1}^k f(x_{i_j}) = \\
&=&\sum_{k=0}^{\infty} \frac{1}{k!} \int_{(0, +\infty)^k} \prod_{j=1}^k f(x_j) \* \rho_k(x_1, \ldots, x_k) dx_1 \cdots dx_k \\
&=&
\sum_{k=0}^{\infty} \frac{1}{k!} \left( \int_{(0, +\infty)} f(x)\* \rho(x) dx \right)^k 
=
\exp\left( \int_{(0, +\infty)} f(x)\* \rho(x) dx \right).
\end{eqnarray}
In the equations above, $\rho_k$ denotes the $k$-point correlation function, and $\rho$ denotes the one-point correlation function 
(also known as intensity). It is a characteristic property of a Poisson random point process that the k-point correlation function factorizes
as a product of one-point correlation functions, i.e. $\rho_k(x_1, \ldots, x_k)=\prod_{i=1}^k \rho(x_i). \ $
In the context of Theorems 2 and 3,  test function $f$  has the form $f(x)= (1 +z\*x)^{-1/2} -1. $ When the intensity $\rho \ $ equals
$ \frac{1}{\pi\* x^{3/2}}, \ $ one obtains 
$$\int_{(0, +\infty)} f(x)\* \rho(x) dx = \int_{(0, +\infty)} ((1 +z\*x)^{-1/2} -1) \* \frac{1}{\pi\* x^{3/2}} dx = 
-\frac{2}{\pi} \* \sqrt{z}, $$  which finishes the proof.

The fact, that the intensity $\rho(x)=\frac{1}{\pi\* x^{3/2}}$ of the Poisson random point process diverges at zero and is 
summable at $+\infty, $
means, that the the vast majority of the eigenvalues of the normalized matrix converge to zero in the limit.

{\bf Remark 2}
{\it It should be pointed out, that the results of Theorem 2 and 3 do not  imply that the statistics of the largest eigenvalues of 
a normalized sample covariance matrix with Gaussian entries are 
Poisson in 
the limit of  $n \to \infty$.  Indeed, to prove the Poisson statistics in the limit one has to show that
\begin{equation}
\label{rav}
\lim_{n \to \infty} E \prod_{i=1}^n \bigl(1+f(\tilde{\lambda_i})\bigr) = {\bf E} \prod_{i=1}^{+\infty}\bigl(1+f(x_i)\bigr) 
\end{equation}
for a sufficiently large class of the test functions $f$, e.g. for step functions with compact support. 
As we already pointed out, the results of Theorems 2 and 3 claim that
(\ref{rav}) is valid for $f(x)= (1+z\*x)^{-1/2} -1 $ for all $z$ such that $ \Re z >0.$}

\section{Conclusion}

It is known in the theory of random Schr\"{o}dinger operators, that the statistics of the eigenvalues is Poisson in the localization regime
(see e.g. \cite{Mo}, \cite{Mi}). It seems, that the same mechanism is responsible for the Poisson statistics for the largest eigenvalues 
in the random matrix models described above. The interesting next problem  is to find a phase transition between the Tracy-Widom 
regime (when all moments of matrix entries exist) and the Poisson regime (when second moment does not exist). 

It is also worth to point out, that there is a vast literature on the Poisson statistics of the energy levels of quantum sysytems in the case 
of the regular underlying dynamics (see e.g. \cite{BT}, \cite{Si}, \cite{Sar}, \cite{CLM}, \cite{Mar1}, \cite{Mar3}).

\printindex
\end{document}